\documentclass[proc]{edpsmath}

\usepackage[T1]{fontenc}

\usepackage{lipsum}
\usepackage{amsfonts, amssymb, dsfont}
\usepackage{graphicx}
\usepackage{epstopdf}
\usepackage{algorithmic}
\usepackage{color}
\usepackage{mathrsfs}
\usepackage[colorlinks = true,
            linkcolor = blue,
            urlcolor  = blue,
            citecolor = blue,
            anchorcolor = blue]{hyperref}
\ifpdf
  \DeclareGraphicsExtensions{.eps,.pdf,.png,.jpg}
\else
  \DeclareGraphicsExtensions{.eps}
\fi

\usepackage{amsopn}
\usepackage{mathtools}
\usepackage{subfigure}
\usepackage[left=0.9in,right=0.9in,top=1in,bottom=1in,includefoot,includehead,headheight=13.6pt]{geometry}

\usepackage[breakable]{tcolorbox}
\newtcolorbox{boxblue}[1]{breakable,colback=white,colframe=BlueMAP,fonttitle=\bfseries,title=#1,left=4pt,right=4pt,top=4pt,bottom=6pt}
\newtcolorbox{boxgreen}[1]{breakable,colback=white,colframe=red,fonttitle=\bfseries,title=#1,left=4pt,right=4pt,top=4pt,bottom=6pt}
\newtcolorbox{partie}{colback=white,colframe=BlueMAP,fonttitle=\bfseries}

\newcommand{\cC}{\ensuremath{\mathcal{C}}}

\newcommand{\cJ}{\ensuremath{\mathcal{J}}}
\newcommand{\cK}{\ensuremath{\mathcal{K}}}

\newcommand{\cR}{\ensuremath{\mathcal{R}}}
\newcommand{\cS}{\ensuremath{\mathcal{S}}}

\newcommand{\cV}{\ensuremath{\mathcal{V}}}

\newcommand{\cX}{\ensuremath{\mathcal{X}}}
\newcommand{\cY}{\ensuremath{\mathcal{Y}}}

\newcommand{\bN}{\ensuremath{\mathbb{N}}}

\newcommand{\bP}{\ensuremath{\mathbb{P}}}

\newcommand{\bR}{\ensuremath{\mathbb{R}}}

\def\[{\left[}
\def\]{\right]}
\def\<{\langle}
\def\>{\rangle}
\def\({\left(}
\def\){\right)}
\def\[{\left [}
\def\]{\right]}
\def\({\left(}
\def\){\right)}

\newcommand{\rd}{\ensuremath{\mathrm d}}
\newcommand{\dx}{\ensuremath{\mathrm dx}}

\newcommand{\eps}{\ensuremath{\varepsilon}}
\newcommand{\e}{\ensuremath{\varepsilon}}

\newcommand{\cond}{\; :\;}

\usepackage{soul}

\begin{document}

\newcommand\unn{u_{\mathcal{NN}}}
\newcommand\znn{z_{\mathcal{NN}}}
\title{Deep Learning-based Schemes \\for Singularly Perturbed Convection-Diffusion Problems}\thanks{This project was partly funded by the Emergence project grant of the Paris City Council ``Models and Measures''. Funding from I-Site FUTURE is also acknowledged. This publication is part of a project that has received funding from
the European Research Council (ERC) under the European Union’s Horizon 2020 Research and
Innovation Programme (Grant Agreement n◦ 810367).}
\author{A.~Beguinet}\address{Paris-Saclay University.}
\author{V.~Ehrlacher}\address{Ecole des Ponts ParisTech \& MATHERIALS INRIA team-project, Marne-la-Vall\'ee, France}
\author{R.~Flenghi}\address{Cermics, \'Ecole des Ponts, INRIA, Marne-la-Vall\'ee, France}
\author{M.~Fuente}\address{COMMEDIA, INRIA Paris, France \& LJLL, Sorbonne Université.}
\author{O.~Mula}\address{Department of Mathematics and Computer Science, Eindhoven University of Technology, 5600 MB, Eindhoven, Netherlands.}
\author{A.~Somacal}\address{LJLL, Sorbonne Université,  Paris, France}
%
%
\begin{abstract}
Deep learning-based numerical schemes such as Physically Informed Neural Networks (PINNs) have recently emerged as an alternative to classical numerical schemes for solving Partial Differential Equations (PDEs). They are very appealing at first sight because implementing vanilla versions of PINNs based on strong residual forms is easy, and neural networks offer very high approximation capabilities. However, when the PDE solutions are low regular, an expert insight is required to build deep learning formulations that do not incur in variational crimes. Optimization solvers are also significantly challenged, and can potentially spoil the final quality of the approximated solution due to the convergence to bad local minima, and bad generalization capabilities. In this paper, we present an exhaustive numerical study of the merits and limitations of these schemes when solutions exhibit low-regularity, and compare performance with respect to more benign cases when solutions are very smooth. As a support for our study, we consider singularly perturbed convection-diffusion problems where the regularity of solutions typically degrades as certain multiscale parameters go to zero.
\end{abstract}
%
%
%
\maketitle





\section{Introduction}
\subsection{Scientific context and goals}
Singularly perturbed differential equations are typically characterized by a small parameter $\e>0$ multiplying some of the highest order terms in the differential equation. In general, the solutions to such equations exhibit multiscale phenomena, and this raises significant challenges to classical numerical methods such as finite elements or finite volumes. To build accurate and robust approximations with these methods as $\e$ decreases, it is necessary  to develop elaborate numerical discretizations. In addition to the mathematical difficulties of the formulation, the resulting numerical schemes are often not entirely trivial to implement: they often require mesh adaptation, and working on complicated geometries is challenging. These difficulties motivate the search for new discretization schemes, hopefully mesh-free, with potential to deliver good quality approximations with easier implementation techniques. In this work, we explore this research direction, and consider strategies based on deep learning techniques. Our main goal is to test various neural network-based schemes, so as to design a strategy which should be robust when $\e\to 0$, easily implementable even for complicated geometries, and with potential to scale in high dimension.

The idea of working with neural network functions to solve PDEs is by far not novel, and countless contributions have been proposed on this front in recent years. The strategies can be roughly classified into two categories:
\begin{enumerate}
\item In the first category, deep neural networks are employed to assist classical numerical methods by improving some limitations, or accelerating certain steps (see, e.g., \cite{LK1990, WM1990, DJ2020}). 
\item In the second category, neural networks are used to directly approximate the solution of PDEs. The solution schemes become in this case an optimization problem where it is crucial to design appropriate loss functions. The loss functions are mostly based on residuals of the equations, and yield to different methods depending on the specific choice:
\begin{enumerate}
\item Physics-informed neural networks (PINNs, \cite{RPK2019}) is a collocation-based method. One finds the coefficients of the neural network solution by minimizing a discretized version of the $L^2$ norm of the strong form of the residual of the PDE. This method is very easily implementable but it implicitly assumes that solutions are very regular.
\item Other strategies leverage weak variational formulations where less regular solutions are allowed. On this front, most of the classical methods originally formulated for piecewise polynomial functions have by now been tested with trial and test spaces of neural network functions.  In this respect, the deep Galerkin method (DGM, \cite{SS2018}) is based on a least-squares formulation, and the variational PINNS (VPINNs, \cite{KZK2019, KZK2021}) is based on the Galerkin method. The main drawback of this approach is that the approximation quality depends on the architecture of both the trial \emph{and} the test neural network classes. In addition, numerous evaluations for multiple test functions need to be performed. Also, strategies involving the minimization of weak-form residuals are usually not trivial to implement because they involve the computation of norms in very weak spaces which necessitate   extra discretization steps.
\item Another approach based on weak variational formulations is the so-called deep Ritz method (DRM, \cite{EY2018}). It leverages the fact that the solution of certain PDEs are the unique minimizer to a certain energy functional. When possible, this approach seems the most appealing: the loss functions is naturally given by the problem, it can accommodate low regular solutions, and the computational cost is moderate in the sense that it only requires to handle test functions (no trial functions). It also carries potential to address high dimensional problem as illustrated in \cite{EHJ2017, EY2018, RMV2020}.
\end{enumerate}
\end{enumerate}

\subsection{Contribution}
The goal of this work is to compare and develop several neural network schemes for singularly perturbed problems when $\e\to 0$. We focus more particularly on convection-diffusion (or stationnary Fokker-Planck) problems with vanishing diffusion for which we explore schemes from the second category according to the above distinction. In other words, we approximate solutions of singular PDEs with feedforward neural network functions. When $\e\to 0$, the regularity of the solutions is deteriorated because of local or boundary thin layers. Therefore the vanilla PINNs method is expected to perform poorly for small values of $\e$ because it commits a variational crime  (and this is actually confirmed in our numerical experiments). Methods based on weak variational formulations seem better adapted, and on that front, it is desirable to work with the deep Ritz method. However, finding energy formulations is not straightforward due to the non-symmetric nature of convective effects. We show how this method can be applied in this context thanks to a change of variable. We compare its numerical robustness with respect to the PINNs method, and a naive finite element discretization with a uniform grid. In the present study, our tests are performed on a 1D example. Despite its simplicity, the example exhibits all the features that are challenging for numerical schemes. For our purposes, the example also presents the important advantage of having analytic solutions which we can leverage in our error analysis, and our validations. Higher-dimensional tests involving also more elaborate sampling strategies are left for future work.

The paper is organized as follows. In Section~\ref{sec:PDE}, various formulations of the convection-diffusion problem we are interested in are introduced. In Section~\ref{sec:neural}, we introduce various neural networks-based schemes which are inspired from the various formulations introduced in Section~\ref{sec:PDE}. The reader is encouraged to observe that an expert mathematical insight is required in order to build formulations that do not incur in variational crimes. Lastly, in Section~\ref{sec:resnum}, these various schemes are compared for one-dimensional problem. We comment on their respective merits and limitations as $\e$  goes to $0$.

\section{A singularly perturbed convection-diffusion equation}\label{sec:PDE}

The aim of this section is to introduce the singularly perturbed convection-diffusion equation we consider in this work, and various formulations of the problem which will be used in Section~\ref{sec:neural} so as to design various neural networks-based schemes for its numerical solution.

\subsection{Problem definition}
As a prototypical example, we consider the following singularly perturbed convection-diffusion equation on a given domain $\Omega\subset \bR^d$, with $d\in \bN^*$. Let $F:\Omega\to \bR^d$ be a given force field, $0< \e \ll 1$ a small parameter, and $f:\Omega\to \bR$ is a given right-hand side function. Our goal is to find a solution $u:\Omega\to \bR$ to
\begin{equation}
\label{general_equation}
-\e (\Delta u)(x) + \nabla \cdot (F u)(x) = f(x), \quad \forall x \in \Omega,
\end{equation}
together with Robin boundary conditions
\begin{equation}\label{Robin}
\alpha (\nabla u \cdot n)(x) + \kappa u(x) = g(x) \quad \forall x \in \partial \Omega,
\end{equation}
where $n$ refers to the outward unit vector of $\partial \Omega$, $\alpha,\,\kappa\geq 0$ and $g$ is a real-valued function defined on $\partial \Omega$. In the following, we assume that the force field $F$ derives from a potential function $V:\Omega\to \bR$, in the sense that
$$
F(x) = - \nabla V(x), \quad \forall x \in \Omega. 
$$
Under appropriate assumptions on $F$ (or $V$), $f$ and $g$, which are assumed to be smooth functions for the sake of simplicity, problem (\ref{general_equation})-(\ref{Robin}) can be proved to have a unique solution~\cite{risken1996fokker,jordan1998variational,bogachev2015fokker}. 
Note that, more generally, $\alpha$ and $\kappa$ could also be given as real-valued functions defined on $\partial \Omega$, instead of constants, and our subsequent developments could be easily adapted.

The equation represents the change in the concentration $u$ of a quantity in a given medium, and in presence of convective and diffusive effects. The force field $F$ represents the drag force while the singular perturbation parameter $\e$ represents the diffusivity of the medium. In the limit of an inviscid medium as $\e\to 0$, the equation changes from elliptic to hyperbolic nature, and from second to first order. For Dirichlet boundary conditions $u=0$ on $\partial \Omega$, the solution can develop sharp boundary layers of width $\e$ near the outflow. We refer the reader to \cite{RST2008} for general references on this equation regarding its analysis and numerical methods.

Classical numerical methods are challenged by problem \eqref{general_equation} when $\e$ is small. In the case of the Galerkin finite element method, the poor performance for this problem is reflected in the bound on the error in the finite element solution. For $\Omega=(0, 1)$ and Dirichlet boundary conditions, a standard Galerkin method with a uniform grid of size $h$ delivers a solution $u_h$ on a finite element space $\bP_h$ that satisfies
\begin{equation}
\Vert u - u_h \Vert_{H^1(0,1)} \leq C(\e) \inf_{w_h \in \bP_h} \Vert u - w_h \Vert_{H^1(0,1)},
\end{equation}
where $C(\e) \sim \e^{-1}$, so that the constant blows up as $\e\to 0$ (see \cite[Theorem 2.49]{RST2008}). The dependence of $C$ on $\e$ is usually referred to as a \emph{loss of robustness} in the sense that, as $\e$ decreases, the Galerkin method is bounded more and more loosely by the best approximation error. As a consequence, on a coarse mesh and for small values $\e$, the Galerkin approximation develops spurious oscillations everywhere in the domain. This very well-known behavior will actually be observed later on in our numerical tests.

Numerous methods have been proposed in order to address this loss of robustness in finite element methods. An important family of methods is based on using residual-based stabilization techniques. Given some variational form, the problem is modified by adding to the bilinear form the strong form of the residual, weighted by a test function and scaled by a stabilization constant $\tau$.  The most well-known example of this technique is the streamline upwind Petrov-Galerkin (SUPG) method (see \cite{BH1982}). The addition of the residual-based stabilization term, can be interpreted as a modification of the test functions which means that these methods seek stabilization by changing the test space, and motivates to search for optimal test spaces in the spirit of \cite{DH2013, CHTD2015}.

Other classical discretization methods such as finite volumes suffer from similar issues, and strategies involving layer-adaptive grids such as Shishkin meshes have been proposed (see, e.g., \cite{KO2010}).

The aim of this work is to explore the potential of approximating solutions of such problems with neural network functions, and the next section presents several options for this, with a discussion on their merits and limitations.

\subsection{General formulation}

Any neural-network based numerical scheme for the solution of (\ref{general_equation})-(\ref{Robin}) relies on the use of a variational formulation of this problem which enables to write $u$ (or another function defined from $u$) as a minimizer of a problem of the form
\begin{equation}\label{genmin}
\mathop{\min}_{v\in \mathcal V} \mathcal J(v), 
\end{equation}
where $\mathcal V$ is a particular set of real-valued functions defined on $\Omega$. The loss function $\cJ:\cV\to \bR$ is usually of the form 
\begin{equation}\label{eq:genJ}
\mathcal J(v):=  \int_\Omega \mathcal R(v)(x) \rd\rho(x) +\int_{\partial \Omega} \mathcal S(v)(r)\rd\tau(r), \quad \forall v\in \cV,
\end{equation}
where for every $v\in \cV$, $\cR(v)$ and $\cS(v)$ are real-valued functions defined on $\Omega$ and $\partial \Omega$ respectively. They are assumed to be measurable with respect to the measures $\rho$ and $\tau$, which are defined on $\Omega$ and $\partial \Omega$ respectively.

The aim of the next sections is to introduce various formulations of problem (\ref{general_equation})-(\ref{Robin}) under the form (\ref{genmin})-(\ref{eq:genJ}). This requires to define appropriate definitions of the set $\mathcal V$, the functions $\mathcal R(v)$ and $\mathcal S(v)$ for any $v\in \mathcal V$ and the unknown function solution of (\ref{genmin}). Unless otherwise stated, the measures $\rho$ and $\tau$ will be defined as the Lebesgue and Lebesgue surfacic measures respectively. 

\subsection{Vanilla (V) formulation}\label{sec:vanilla}

We begin by introducing the most classical formulation used in neural network-based numerical schemes such as PINNs. For the reasons that we outline next, different aspects of this formulation can be improved, therefore we refer to it as \textit{vanilla} (V) formulation in the following. 

The formulation consists in interpreting the solution $u$ of \eqref{general_equation}-\eqref{Robin} as the unique solution of a minimization problem of the form (\ref{genmin}) with $\mathcal V= H^2(\Omega)$ and to define for all $v\in \mathcal V$,
\begin{equation}
\begin{cases}
\mathcal R (v)(x) :=\lambda \left| -\e (\Delta v)(x) + \nabla \cdot (F v)(x) - f(x)\right|^2, & \quad \mbox{ for all }x\in \Omega,\\
\mathcal S (v)(x):=(1-\lambda) \left| \alpha (\nabla v \cdot n)(x) + \kappa v(x)- g(x)\right|^2, & \quad \mbox{ for all }x\in \partial\Omega,
\end{cases}
\label{eq:loss-vanilla}
\end{equation}
for some $\lambda \in (0,1)$. In this approach, the parameter $\lambda$ enables to tune the respective weight of the contributions of the bulk and boundary terms in the total functional $\mathcal J$ to be minimized. In practice, in the numerical tests presented in Section~\ref{sec:resnum}, $\lambda$ will always be chosen to be equal to $\frac{1}{2}$.
 
Note that such an approach requires the solution $u$ to belong to $H^2(\Omega)$, which implies that the solution has to be sufficiently regular. When $\eps\to0$, this assumption becomes less and less realistic due to the formation of boundary layers. This raises the question as to whether it is possible to introduce another formulation of problem (\ref{general_equation})-(\ref{Robin}) which would allow for less regular solutions. The goal of the next section is to introduce such an alternative formulation.

\subsection{Weak variational (W) formulation}
\label{sec:variational}
In this section we develop an avenue based on an energy minimization approach which requires less regularity in the solutions than the vanilla formulation. To this aim, we introduce the change of variable
\begin{equation}
\label{eq:change-var}
u(x) = e^{cV(x)}z(x),
\end{equation}
where $c\in \bR$ is a constant yet to be determined. Taking first and second derivatives in \eqref{eq:change-var} yield that for all $x\in \Omega$,
\begin{align*}
\nabla u(x) &= e^{cV(x)} \left( c\nabla V(x)z(x)+ \nabla z(x) \right) \\
\Delta u(x) &= e^{cV(x)} \left( c\Delta V(x) z(x)+|c\nabla V(x)|^2z(x)+2c\nabla V(x) \cdot \nabla z(x)+\Delta z(x) \right).
\end{align*}
Now, setting the value of $c$ to be
$$
c= \frac{1}{2\e},
$$
and inserting the change of variable into \eqref{general_equation}, we conclude that $u$ is a solution to \eqref{general_equation} if and only if $z$ is a solution to the elliptic problem
\begin{align}
\label{eq:v}
- \Delta z(x)+ \left( \frac{\Delta V(x)}{2\e} + \frac{|\nabla V(x)|^2}{4\e^2}  \right) z(x)
&= f(x)\frac{e^{-\frac{V(x)}{2\e}}}{\e}, \quad \forall x\in \Omega, 
\end{align}
with Robin boundary conditions
\begin{equation}
\label{eq:robin-v}
\alpha (\nabla z(x) \cdot n(x)) + \left( \kappa + \frac{\alpha}{2\e} \nabla V(x) \cdot n(x) \right) z(x) = e^{\frac{-V(x)}{2\e}}g(x),\quad \forall x \in \partial \Omega .
\end{equation}
At this stage, one could of course apply the vanilla formulation to solve \eqref{eq:v}-\eqref{eq:robin-v} and compute $z$ solution of a minimization problem of the form (\ref{genmin}) with $\cV = H^2(\Omega)$ and the functionals $\mathcal R$ and $\mathcal S$ defined by
\begin{equation}
\begin{cases}
\mathcal R (v)(x) :=  \lambda\left|\Delta v(x)+ \left( \frac{\Delta V(x)}{2\e} + \frac{|\nabla V(x)|^2}{4\e^2}  \right) v(x)-f(x)\frac{e^{-\frac{V(x)}{2\e}}}{\e}\right|^2, & \quad \mbox{ for all }x\in \Omega,\\
\mathcal S (v)(x):=\left|(1-\lambda) \alpha (\nabla v(x) \cdot n(x)) + \left( \kappa + \frac{\alpha}{2\e} \nabla V(x) \cdot n(x) \right) v(x)- e^{\frac{-V(x)}{2\e}}g(x)\right|^2, & \quad \mbox{ for all }x\in \partial\Omega,
\end{cases}
\label{eq:loss-cv-vanilla}
\end{equation}
for all $v\in \cV = H^2(\Omega)$ and some $\lambda \in (0,1)$. The value of $\lambda$ chosen in our numerical tests is $\lambda =0.5$. We will refer to this approach as the \itshape vanilla-z \normalfont ($Vz$) formulation.   

    \medskip

Note that this method does not fully exploit the change of variables since the elliptic nature of problem \eqref{eq:v} allows us to easily build a weak formulation of this equation. Testing against a smooth test function $v$ and integrating by parts we obtain the weak formulation
\begin{align*}
&\int_{\Omega} \nabla z \cdot \nabla v
- \int_{\partial \Omega} v \nabla z \cdot n \dx
+ \int_\Omega \left( \frac{\Delta V}{2\e} + \frac{|\nabla V|^2}{4\e^2}  \right) z v   = \int_\Omega  f\frac{e^{-\frac{V}{2\e}}}{\e} v.
\end{align*}
Using equality \eqref{eq:robin-v}, we get
\begin{align*}
\int_{\Omega} \nabla z \cdot \nabla v
+ \int_\Omega \left( \frac{\Delta V}{2\e} + \frac{|\nabla V|^2}{4\e^2}  \right) z v  
+ \int_{\partial \Omega} \left( \frac{\kappa}{\alpha}+ \frac{1}{2\e} \nabla V \cdot n \right) zv\\
= \int_\Omega  f\frac{e^{-\frac{V}{2\e}}}{\e} v  +  \int_{ \partial \Omega} \frac{1}{\alpha} e^{\frac{-V}{2\e}}g v .
\end{align*}
Therefore the weak formulation of problem \eqref{eq:v} is to find $z\in H^1(\Omega)$ such that
\begin{equation}
\label{eq:bilin-v}
a(z, v) = \ell (v),\quad \forall v\in H^1(\Omega)
\end{equation}
with
\begin{align*}
a(z, v) &\coloneqq \int_{\Omega} \nabla z \cdot \nabla v
+ \int_\Omega \left( \frac{\Delta V}{2\e} + \frac{|\nabla V|^2}{4\e^2}  \right) zv \dx  
+  \int_{\partial \Omega} \left(\frac{\kappa}{\alpha}+ \frac{1}{2\e} \nabla V \cdot n \right) zv
\\
l(v) &\coloneqq \int_\Omega  f\frac{e^{-\frac{V}{2\e}}}{\e} v +  \int_{\partial \Omega}  \frac{1}{\alpha} e^{\frac{-V}{2\e}}gv 
\end{align*}
To ensure that the symmetric bilinear form $a$ is continuous and coercive, we assume in the sequel that the following conditions are satisfied:
\begin{align}
\label{eq:cond-v}
\begin{cases}
\left( \frac{\Delta V(x)}{2\e} + \frac{|\nabla V(x)|^2}{4\e^2}  \right) \geq a_0 > 0, &\quad \forall x \in \Omega \\
\left( \frac{\kappa}{\alpha}+ \frac{1}{2\e} \nabla V(x) \cdot n(x) \right) \geq 0, &\quad\forall x \in \partial \Omega
\\
f\in L^2(\Omega),\;g\in L^2(\partial\Omega)
\end{cases}
\end{align}
In that case, all the hypothesis of the Lax-Milgram theorem are satisfied for $a$ and $\ell$. Thus, $z$ can be equivalently rewritten as the unique solution of a minimization problem of the form
\begin{equation}\label{eq:var}
z = \mathop{\rm argmin}_{v\in H^1(\Omega)} \frac{1}{2} a(v,v) - \ell(v). 
\end{equation}

This implies that $z$ can equivalently be recast as the unique solution of a minimization problem of the form
(\ref{genmin}) with $\cV = H^1(\Omega)$ and 
\begin{equation}
\label{eq:loss-variational-detail}
\begin{cases}
\mathcal R(v)(x) &\coloneqq \frac{1}{2} \left[ |\nabla v(x)|^2 + \left( \frac{\Delta V(x)}{2\e} + \frac{|\nabla V(x)|^2}{4\e^2}  \right) |v(x)|^2\right] - f(x)\frac{e^{-\frac{V(x)}{2\e}}}{\e} v(x), \quad \forall x\in \Omega, 
\\
\mathcal S(v)(x)&\coloneqq \frac{1}{2}\left[ \left(\frac{\kappa}{\alpha}+ \frac{1}{2\e} \nabla V(x) \cdot n(x) \right) |v(x)|^2 \right] - \frac{1}{\alpha} e^{\frac{-V(x)}{2\e}}g(x) v(x) , \quad \forall x\in \partial \Omega. 
\end{cases}
\end{equation}
We will refer to this approach as the \itshape weak-z \normalfont (Wz) formulation.

\medskip

Moreover, using (\ref{eq:change-var}), we can equivalently express $u$ as a solution of a minimization problem of the form (\ref{genmin}) with 
\begin{equation}\label{eq:H1CV}
\cV:= \left\{ v= e^{\frac{V}{2\e}} \overline{v}, \;  \overline{v}\in H^1(\Omega)\right\}, 
\end{equation}
and 
\begin{equation}
\label{eq:loss-variational-detail2}
\begin{cases}
\mathcal R(v)(x) &\coloneqq \frac{1}{2} \left[ |\nabla \overline{v}(x)|^2 + \left( \frac{\Delta V(x)}{2\e} + \frac{|\nabla V(x)|^2}{4\e^2}  \right) |\overline{v}(x)|^2\right] - f(x)\frac{e^{-\frac{V(x)}{2\e}}}{\e} \overline{v}(x), \quad \forall x\in \Omega, 
\\
\mathcal S(v)(x)&\coloneqq \frac{1}{2}\left[ \left(\frac{\kappa}{\alpha}+ \frac{1}{2\e} \nabla V(x) \cdot n(x) \right) |\overline{v}(x)|^2 \right] - \frac{1}{\alpha} e^{\frac{-V(x)}{2\e}}g(x) \overline{v}(x) , \quad \forall x\in \partial \Omega, 
\end{cases}
\end{equation}
for all $v\in \mathcal V$ with $\overline{v}:= v e^{-\frac{V}{2\e}}$. We will refer to this formulation as the \itshape weak \normalfont (W) formulation.  

\subsection{Rescaled formulation}

In this section, we introduce another formulation based on a change of scale in the original problem. More precisely, introducing $\Omega_\e:= \frac{1}{\e}\Omega$, we introduce auxiliary functions $\widetilde{u}: \Omega_\e \to \mathbb{R}$, $\widetilde{z}: \Omega_\e \to \mathbb{R}$ and $\widetilde{V}: \Omega_\e \to \mathbb{R}$ defined so that for all $x\in \Omega$, 
\begin{equation}\label{rescaling}
u(x) =\e \widetilde{u}\left(\frac{x}{\e}\right),\quad z(x) =\e \widetilde{z}\left(\frac{x}{\e}\right),\quad V(x) =\e \widetilde{V}\left(\frac{x}{\e}\right).
\end{equation}

Notice that if $u$ and $z$ satisfy \eqref{eq:change-var}, then
$$
\widetilde{z}(y) = \widetilde{u}(y) e^{\frac{1}{2}\widetilde{V}(y)}, \quad \forall y \in \Omega_\e. 
$$
Denoting by $\widetilde{F}(y):= - \nabla \widetilde{V}(y) = F\left(\e y\right)$ for all $y\in \Omega_\e$, it holds that $u$ is solution to (\ref{general_equation})-(\ref{Robin}) if and only if $\widetilde{u}$ is solution to 
\begin{equation}\label{PDErescaled}
- \Delta \widetilde{u}(y) +  \nabla \cdot \left( \widetilde{F} \widetilde{u}\right)(y)  = \widetilde{f}(y), \quad \forall y \in \Omega_\e
\end{equation}
where $\widetilde{f}(y):= f(\e y)$ for all $y\in \Omega_\e$ with boundary conditions
\begin{equation}\label{Robinrescaled}
\alpha (\nabla \widetilde{u}\cdot n)(y) + \e \kappa \widetilde{u}(y) = \widetilde{g}(y), \quad \forall y\in \partial \Omega_\e, 
\end{equation}
with $\widetilde{g}(y):= g(\e y)$ for all $y\in \Omega_\e$.

Using similar calculations to the ones done in Section~\ref{sec:variational}, the Lax-Milgram theorem guarantees that $\widetilde{z}$ is the unique solution in $H^1(\Omega_\e)$ of the following variational problem: for all $\widetilde{v} \in H^1(\Omega_\e)$,
\begin{align}
&\int_{\Omega_\e} \nabla \widetilde{z}(y) \cdot \nabla \widetilde{v} (y) \,dy
+ \int_{\partial \Omega_\e} \left(\frac{\e\kappa}{\alpha}+\frac{1}{2}\nabla\tilde{V}(y)\cdot n(y)\right) \widetilde{z}(y) \widetilde{v}(y) \,dy
+ \int_{\Omega_\e} \left(\frac{\Delta \tilde{V}\left(y\right)}{2}+\frac{|\nabla \tilde{V}\left(y\right)|^2}{4}\right)\widetilde{z}\left(y\right) \widetilde{v}(y) \,dy \nonumber\\&
\qquad= \int_{\Omega_\e}  \widetilde{f}(y)e^{-\frac{1}{2}\tilde{V}(y)}\widetilde{v}(y) \,dy+\int_{\partial\Omega_\e}\frac{1}{\alpha}e^{-\frac{1}{2}\tilde{V}\left(y\right)}\widetilde{g}(y)\widetilde{v}(y)\,dy.
\label{eq:var-rescaled}
\end{align}
The result is valid provided that the following assumptions on the coefficients hold
\begin{align*}
\begin{cases}
\Delta \tilde{V}(y)+|\nabla\tilde{V}|^2(y)\geq 0, \; \forall y \in \Omega_\e,\\
  \frac{\kappa \e}{\alpha}+\frac{1}{2}\nabla\tilde{V}(y)\cdot n(y)\geq 0, \; \forall y \in \partial \Omega_\e,\\
  \widetilde{f}\in L^2(\Omega_\e),\;\widetilde{g}\in L^2(\partial\Omega_\e).
 \end{cases}
\end{align*}
The above conditions are equivalent to the assumptions \eqref{eq:cond-v} stated in Section~\ref{sec:variational}.

\medskip

As in the previous section, the function $\widetilde{z}$ is then the unique solution of a minimization problem of the form (\ref{genmin}) with $\mathcal V = H^1(\Omega_\e)$ with 
\begin{equation}
\label{eq:rescaled}
\begin{cases}
\mathcal R(v)(y) &\coloneqq \frac{1}{2} \left[ |\nabla v(y)|^2 + \left( \frac{\Delta \widetilde{V}(y)}{2} + \frac{|\nabla \widetilde{V}(y)|^2}{4}  \right) |v(y)|^2\right] - \widetilde{f}(y)e^{-\frac{\widetilde{V}(y)}{2}} v(y), \quad \forall y\in \Omega_\e, 
\\
\mathcal S(v)(y)&\coloneqq \frac{1}{2}\left[ \left(\frac{\e\kappa}{\alpha}+ \frac{1}{2} \nabla \widetilde{V}(y) \cdot n(y) \right) |v(y)|^2 \right] - \frac{1}{\alpha} e^{\frac{-\widetilde{V}(y)}{2}}\widetilde{g}(y) v(y) , \quad \forall y\in \partial \Omega_\e. 
\end{cases}
\end{equation}

We will refer to this approach as the \textit{rescaled-weak-z} \normalfont (RWz) formulation.

\subsection{Summary of the methods}
For the sake of clarity, we summarize here the main features of each method.

\begin{center}
\begin{tabular}{|c||c|c|c|c|}
\hline
 Method & Acronym & Unknown &  $\mathcal V$ & $\mathcal R$ and $\mathcal S$\\
 \hline\hline
 vanilla & V & $u$ & $H^2(\Omega)$ & \eqref{eq:loss-vanilla}\\
 \hline
 vanilla-z & Vz & $z$ & $H^2(\Omega)$ & \eqref{eq:loss-cv-vanilla}  \\
 \hline
 weak & W & $u$ & (\ref{eq:H1CV}) & \eqref{eq:loss-variational-detail2}\\
 \hline
  weak-z & Wz & $z$ & $H^1(\Omega)$ & \eqref{eq:loss-variational-detail}\\
 \hline
 rescaled-weak-z & RWz & $\widetilde{z}$ & $H^1(\Omega_\eps)$ & \eqref{eq:rescaled} \\
 \hline
\end{tabular}
\end{center}

\section{Neural networks based numerical schemes}\label{sec:neural}

In this section we describe the numerical approach used in order to compute an approximation of the solution of a minimization problem of the form (\ref{genmin}) by means of a neural-network based method. We first present in Section~\ref{sec:general} the general principle of such approaches.  The main ingredients to design a neural-network based method consist in the choice of a class of neural network functions and of sampling schemes in order to approximate the integrals involved in the definition of the loss function $\mathcal J$ defined by (\ref{eq:genJ}). These two ingredients are detailed respectively in Section~\ref{sec:neural-nets} and Section~\ref{sec:sampling} respectively. Finally, some details on the numerical implementation are given in Section~\ref{sec:details}.

\subsection{General principle}\label{sec:general}

The numerical solution of a minimization problem of the form \eqref{genmin} usually requires to consider alternatives to $\cV$ and $\cJ$ that are amenable for practical implementation. The strategy thus consists in formulating a related problem of the form

\begin{equation}\label{minapp}
\mathop{\min}_{v\in \mathcal K} \widehat{\mathcal J}(v), 
\end{equation}
where 
\begin{itemize}
    \item $\cK \subset \cV$ is a set of functions parametrized by a finite number of scalar coefficients. A classical class of functions are finite elements. Here, we consider neural networks (see Section~\ref{sec:neural-nets} below); 
    \item $\widehat{\mathcal J}$ is an approximation of the loss function $\mathcal J$ where the integrals are approximated using some particular quadrature or sampling schemes. 
\end{itemize}

More precisely, for given integers $K, M \in \mathbb{N}^*$, given sets of points $(x_k)_{1\leq k \leq K} \subset \Omega$, $(y_m)_{1\leq m \leq M} \subset \partial \Omega$, and given sets of weights $(\rho_k)_{1\leq k \leq K}\subset \mathbb{R}_+$ and $(\tau_m)_{1\leq m \leq M}\subset \mathbb{R}_+$, for all $v\in \mathcal K$, the functional $\widehat{\mathcal J}(v)$ is defined by
\begin{equation}\label{eq:Jhatdef}
\widehat{\mathcal J}(v):= \sum_{k=1}^K \rho_k \mathcal R(v)(x_k) + \sum_{m=1}^M \tau_m \mathcal S(v)(y_m). 
\end{equation}

As a consequence, the definition of a neural-network based numerical scheme for the approximation of a problem of the form (\ref{genmin}) requires the definition of two ingredients: 
\begin{itemize}
    \item the class $\mathcal K \subset \mathcal V$ of neural network functions; 
    \item the sampling scheme, i.e. the choice of $K$, $M$, $(x_k)_{1\leq k \leq K}$, $(y_m)_{1\leq m \leq M}$, $(\rho_k)_{1\leq k \leq K}$ and $(\tau_m)_{1\leq m \leq M}$ in order to define the approximate functional $\widehat{\mathcal J}$ given by (\ref{eq:Jhatdef}). 
\end{itemize}

The set of neural network functions $\mathcal K$ we consider in our numerical experiments is presented in Section~\ref{sec:neural-nets}. The various sampling schemes tested here are given in Section~\ref{sec:sampling}.



\subsection{Neural Network classes of functions}
\label{sec:neural-nets}

In this work, we only consider classes of functions defined by means of feedforward neural networks whose definition we recall next (see~\cite{ma2020towards} for general references).

Let $\cX\subset \mathbb{R}^{d_{\mathcal X}}$ and $\cY\subset \bR^{d_{\cY}}$ be some input and output sets of finite dimensions  $d_\cX, d_\cY \in \mathbb{N}^*$. A feedforward neural network is a function
$$
\psi : \cX \to \cY
$$
which reads as
\begin{equation}\label{eq:defpsi}
\psi (x) = T_{L}(\sigma(T_{L-1}(\sigma(\ldots \sigma(T_0(x))))), \quad \forall x \in \cX.
\end{equation}
For every $\ell \in \{0, \dots, L\}$,
\begin{equation}
T_\ell : \left\{ \begin{array}{ccc}
\bR^{p_\ell} &\to & \bR^{p_{\ell+1}} \\
x_\ell &\mapsto & T_\ell(x_\ell) \coloneqq A_\ell x_\ell + b_\ell\\
\end{array}
\right.
\end{equation}
is an affine function which can be expressed through a matrix $A_\ell \in \bR^{p_{\ell+1}\times p_\ell}$, and an offset vector $b_\ell \in \bR^{p_{\ell+1}}$, and $\sigma: \bR \to \bR$ is called the (nonlinear) activation function. By a slight abuse of notation, for all $p\in \mathbb{N}^*$ and for any vector $w:=(w_i)_{1\leq i \leq p}\subset \mathbb{R}^p$, the notation $\sigma(w)$ actually denotes the vector of $\bR^p$ with entries $\sigma(w_i)$, that is, $\sigma(w) = (\sigma(w_i))_{i=1}^p$. Note that since $\psi$ maps $\cX$ onto $\cY$, it is necessary that $p_0 = d_\cX$ and $p_{L+1}=d_\cY$. The layers numbered from $1$ to $L$ are usually called the hidden layers of the neural network.

To define a class of feedforward neural networks, we fix an architecture by prescribing a given activation function $\sigma$, depth $L\in \bN$, and layer widths  $\boldsymbol{p} =(p_0, \dots, p_{L+1})\in (\mathbb{N}^*)^{L+2}$. Once the values of $\sigma$, $L$ and $\boldsymbol{p}$ have been chosen, we view the coefficients $(A_\ell, b_\ell)_{0\leq \ell \leq L}$ of the affine mappings $T_0, \cdots, T_L$ as parameters. We gather these coefficients in the vector of parameters
$$
\theta \coloneqq \{(A_\ell, b_\ell)\}_{\ell=0}^L,
$$
and assume that $\theta$ takes values in a set 
$$\Theta\subseteq \bigtimes_{\ell = 0}^L \left(\mathbb{R}^{p_\ell \times p_{\ell+1}} \times \mathbb{R}^{p_{\ell+1}}\right).
$$
For any $\theta \in \Theta$, we define by $\psi_\theta: \mathcal X\to \mathcal Y$ the function $\psi$ defined by (\ref{eq:defpsi}) with $\theta = \{(A_\ell, b_\ell)\}_{\ell=0}^L\in \Theta$. 

The class of neural network functions with architecture $(\sigma, L, \boldsymbol{p})$ and coefficient sets $\Theta$ is then defined as
$$
\mathcal N(\sigma, L, \boldsymbol{p}, \Theta):= \left\{ \psi_\theta\cond \theta \in \Theta\right\}. 
$$
In our context, the input and output sets $\mathcal X$ and $\mathcal Y$ are respectively given by 
$$
\mathcal X = \Omega \mbox{ (or }\Omega_\e \mbox{)} \quad \mbox{ and } \quad \mathcal Y = \mathbb{R},
$$
so that $d_{\mathcal X} = d$ and $d_{\mathcal Y} = 1$. In all the numerical tests presented below, the class $\mathcal K$ is chosen as 
$$
\mathcal K \coloneqq  \mathcal N(\sigma, L, \boldsymbol{p}, \Theta),
$$
with
$$
\sigma = \tanh,\; L=2, \;  \boldsymbol{p}=(10, 10) \text{ and }\Theta = \bigtimes_{\ell = 0}^L \left(\mathbb{R}^{p_\ell \times p_{\ell+1}} \times \mathbb{R}^{p_{\ell+1}}\right).
$$

Note that the set $\mathcal K$ is then a subset of $\mathcal V$ for all the formulations of the convection-diffusion problem we introduced in Section~\ref{sec:PDE}. Moreover, the solution of the approximate problem (\ref{minapp}) is equivalent to finding a minimizer $\theta^*\in \Theta$ solution to 
\begin{equation}\label{eq:thetastar}
\min_{\theta \in \Theta} \widehat{\mathcal J}(\psi_\theta).
\end{equation}

\textbf{Remark:} In many machine learning applications, the choice of relu activation functions is very common due to its low computational cost when performing evaluation or first order differentiation. However, in our problem, second order derivatives are needed to calculate the loss function. If relu activation functions were used, then the second order derivative terms would be $0$, and no good approximation could be learned. This reason motivates our choice of $\tanh$ as the activation function.

\subsection{Sampling schemes}\label{sec:sampling}

We detail in this section the various sampling schemes we considered in our numerical tests in order to define the approximate loss function $\widehat{\cJ}$.

Since we work with one-dimensional examples, we carry the discussion for dimension one. In fact, we consider problem (\ref{general_equation})-(\ref{Robin}) with $\Omega = (0,1)$ so that $\partial \Omega = \{0\} \cup \{1\}$ (and $\partial \Omega_\e = \{0\} \cup \{1/\e\}$). Thus, for all our tests, the domain boundary has $M = 2$ points $y_1 = 0$ and $y_2 = 1$ (or $y_2 = \frac{1}{\e}$ for the RWz method). Taking $\tau_1 = \tau_2=1$ for the surface weights, the surface term in \eqref{eq:Jhatdef} takes the simple form
$$
\sum_{m=1}^{M=2} \tau_m \mathcal S(v)(y_m) = \int_{\partial \Omega} \mathcal S(v)\,d\tau \quad \forall v \in \cV \mbox{  (or }\int_{\partial \Omega_\e} \mathcal S(v)\,d\tau\mbox{ for the RWz formulation).}
$$

We consider three different sampling schemes for the approximation of the bulk term $\int_\Omega \mathcal R(v)\,d\rho$:
\begin{enumerate}
\item The first choice is a simple \itshape uniform \normalfont sampling scheme (labeled $-u$ in our tests). For a given $K\in \mathbb{N}^*$, we set $\rho_k = \frac{1}{K}$ and $(x_k)_{1\leq k \leq K}$ as the centers of the intervals given by a uniform discretization grid of the interval $(0,1)$.
\item The second sampling scheme, called \itshape random \normalfont ($-r$) scheme, consists in choosing $\rho_k = \frac{1}{K}$ and the points $(x_k)_{1\leq k \leq K}$ as a collection of random points, identically independently distributed according to the uniform distribution on $(0,1)$. 
\item We lastly consider a third sampling scheme, called \itshape exponential \normalfont (-e) scheme, which is specific to the Wz formulation. Recall that in this case, for all $v\in \mathcal K$, the expression of $\mathcal R(v)$ is given by (\ref{eq:loss-variational-detail}), namely
\begin{equation*}
\cR (v)(x) =  \cR^{(1)} (v)(x) + \cR^{(2)} (v)(x), \quad 
\forall x\in \Omega
\end{equation*}
with
\begin{equation*}
\begin{cases}
\cR^{(1)} (v)(x) &\coloneqq \frac{1}{2} \left[ |\nabla v(x)|^2 + \left( \frac{\Delta V(x)}{2\e} + \frac{|\nabla V(x)|^2}{4\e^2}  \right) |v(x)|^2\right] \\
 \cR^{(2)} (v)(x) &\coloneqq - f(x)\frac{e^{-\frac{V(x)}{2\e}}}{\e} v(x), \quad \forall x\in \Omega.
\end{cases}
\end{equation*}
The thus view the bulk integral term as
$$
\int_\Omega \mathcal R(v)(x) \rd \rho(x)
=
\int_\Omega \mathcal R^{(1)}(v)(x) \rd \rho(x)
+
\int_\Omega \mathcal R^{(2)}(v)(x) \rd \rho(x),
$$
and we approximate each component separately as follows. For the first term, we draw a collection of $K_1\in \mathbb{N}^*$ iid.~random points $(x^{(1)}_k)_{1\leq k \leq K_1}$ from the uniform distribution on $(0,1)$ and for all $1\leq k \leq K_1$, the weights $\rho_k^{(1)}$ are chosen to be equal to $\frac{1}{K_1}$. For the second term, we draw $K_2\in \mathbb{N}^*$ iid random points $(x^{(2)}_k)_{1\leq k \leq K_2}$ following the probability density
$$
\rho^{(2)}(x):= \frac{e^{-\frac{V(x)}{2\e}}}{Z_\e}, \quad x\in \Omega,
$$
with 
$$
Z_\e:= \int_\Omega e^{-\frac{V(x)}{2\e}}\,dx
$$
Setting now $\rho_k^{(2)} = \frac{Z_\e}{K_2}$ for all $1\leq k \leq K_2$, the integral $\int_\Omega \mathcal R(v)$ is the approximated by
$$
\sum_{k=1}^{K_1} \rho_k^{(1)} \left( \frac{1}{2} \left[ |\nabla v(x_k^{(1)})|^2 + \left( \frac{\Delta V(x_k^1)}{2\e} + \frac{|\nabla V(x_k^{(1)})|^2}{4\e^2}  \right) |v(x_k^{(1)})|^2\right]\right) - \sum_{k=1}^{K_2} \rho_k^{(2)} \left( f(x_k^{(2)})\frac{1}{\e} v(x_k^{(2)}) \right).
$$
\end{enumerate}

In the following, we use the notation $-u$ (respectively $-r$ and $-e$), after the name of a formulation, in order to refer to the numerical method obtained by using this formulation, together with a uniform (respectively random or exponential) sampling scheme. For instance, the $V-u$ method refers to the vanilla formulation used in conjunction with a uniform sampling scheme.

\subsection{Comparison with finite element schemes}\label{sec:FEM}

One important point in the investigation of the merits and limitations of deep learning-based numerical schemes is to understand how they compare with respect to other existing  schemes. In our tests, we provide a numerical comparison with a vanilla finite element Galerkin scheme involving a uniform mesh. For the sake of completeness, we briefly recall the main steps of our finite element Galerkin approach.

Integrating the original equation \eqref{general_equation} against a sufficiently smooth function $v\in \cC^\infty(\Omega)$, and integrating by parts, it follows that a weak formulation of problem \eqref{general_equation} is to find $u\in H^1(\Omega)$ such that
$$
a(u, v) = l(v), \quad \forall v\in H^1(\Omega)
$$
with
\begin{align}\label{eq:fem u}
a(u,v)&=\int_\Omega\nabla u \nabla v\dx+\e^{-1}\int_\Omega \nabla\cdot(Fu)v\dx+\frac{\kappa}{\alpha}\int_{\partial\Omega}u(x)v(x)\dx \\
l(v)&=\e^{-1}\int_\Omega f(x)v(x)\dx-\int_{\partial\Omega} g(x)v(x)\dx.
\end{align}
We numerically solve this problem by Galerkin projection. For this, we consider a mesh $(T_n)_{n=1}^N$ of $\Omega$ and define the associated $\bP_1$ finite element space
\begin{equation*}
\cV_N\coloneqq\{v\in \mathscr{C}^0( \mathbb{R})\cond\forall 0\leq s\leq N-1,\; v_{[x_s,x_{s+1})}\in\mathbb{P}^1([x_s,x_{s+1}))\}
\end{equation*}
with 
\begin{equation*}
\mathbb{P}^1([x_s,x_{s+1}]\coloneqq \{v:[x_s,\;x_{s+1]},\;v(x)=ax+b,\;(a,b)\in\mathbb{R}^2\}.
\end{equation*}
We then search for a solution $u_N \in \cV_N \subset H^1(\Omega)$ by Galerkin projection, that is, we search for $u_N \in \cV_N$ such that
$$
a(u_N, v) = l(v), \quad \forall v \in \cV_N.
$$
We next take as a basis of $\cV_N$ the set of tent functions defined as 
$$
\varphi_i(x_j)=\delta_{ij},\; \mathrm{for}\; 1\leq i,j\leq N.
$$
and we express the solution as $u_N = \sum_{i=1}^N c_i \varphi_i$. Gathering the expansion coefficients in the vector $c=(c_i)_{i=1}^N$, and injecting the expansion of $u_N$ in the variational formulation, we are led to the system of equations
$$
M c = q
$$
where
\begin{align*}
M &= (M_{i, j})_{1\leq i, j \leq N}, \quad M_{i, j} \coloneqq a(\varphi_i, \varphi_j), \\
q &= (q_i)_{i = 1}^N, \quad q_i \coloneqq l(\varphi_i).
\end{align*}

\section{Numerical Results}\label{sec:resnum}

\subsection{Test case and comparison criteria}
In this section we show the results obtained by approximating the exact solution of the problem described in equation \eqref{general_equation} using the methods introduced above. We work on the one dimensional domain $\Omega=(0,1)$ with $F=1$ and $f=1$. We choose Robin boundary conditions that mimic Dirichlet conditions and we set $\alpha = 10^{-3}$, $\kappa = 1$, $g_0=g_1=0$. Note that we cannot take $\alpha=0$ since all variational methods are not well defined for pure Dirchlet boundary conditions. With these choices, the equation reads 
\begin{equation}\label{eq:dim1}
\left\{
\begin{array}{rl}
- \e u''(x) + u'(x) = 1, & \quad \forall x\in (0, 1),\\
 - 10^{-3} u'(0) +  u(0) = 0, & \\
 10^{-3} u'(1) + u(1) = 0. &
\end{array}\right. 
\end{equation}
Since we work in dimension 1, we can benefit from the fact that the exact solution $u$ has an analytic form as shown in the Appendix \ref{sec:AppendixAnaliticSolution}. We can thus easily compare the approximation quality of the output functions $\hat u$ from our methods by computing a discrete version of their $L^2(\Omega)$ error norm with respect to the exact solution. 
$$
e^2_{L^2} =
\Vert u - \hat u \Vert^2_{L^2(\Omega)}
\approx
\frac{1}{\tilde K}\sum ^{\tilde K-1}_{k=0}(u(x_k) - \hat{u}(x_k))^2
= e^2_{\ell^2}
$$
The points $x_k$ are sampled uniformly as defined in \ref{sec:sampling}. We use $10$ times more points than the ones used for training, so $\tilde K = 10 K$. Similarly, we also compute the error with respect to the $H^1(\Omega)$ semi-norm.
$$
e^2_{H^1} =
\Vert u' - \hat u' \Vert^2_{L^2(\Omega)}
\approx
\frac{1}{\tilde K}\sum ^{\tilde K-1}_{k=0}(u'(x_k) - \hat{u}'(x_k))^2
= e^2_{h^1}
$$
Note that one can obtain the $H^1$ error by adding the above error components.

We study the impact on the errors of the following parameters:
\begin{itemize}
    \item The values of $\e$. They range from $5.10^{-3}$ to $10.0$ with a logarithmic spacing.
    \item The number $K$ of training points. We consider $K=10, \, 10^2, \,10^3,\,10^4$.
    \item The choice of the sampling method for the training points (uniformly spaced or uniformly random, labelled as $-u$ and $-r$).
    \item The impact of the machine precision (Float16, Float32, Float64).
\end{itemize}
Due to the randomness in the initialization of weights on the neural networks, for each combination of parameters ($\e$, $K$, sampling type, and machine precision), we perform 10 repetitions with different initializations. Since we didn't notice a big difference between the $l^2$ error and the $h^1$ error, we keep just the second one for clarity and put in the Appendix \ref{sec:l2error} the plots in $l^2$ error.

\subsection{Our code and practical implementation details}\label{sec:details}

All our neural network based numerical tests were performed in Python 3.6 and using the TensorFlow 1.13.1 library \cite{Tensorflow}. The code provided in the original paper on PINNs \cite{PINNs} was used as the starting point for our own code developments, and we have followed similar guidelines to generalize and enlarge it where needed. In the same way, for each numerical method, derivatives of functions $v\in\mathcal{K}$ are computed using automatic differentiation. The numerical optimization procedure used in order to compute an approximation of $\theta^*$ a minimizer of problem (\ref{eq:thetastar}) is given by the quasi-Newton L-BFGS algorithm \cite{BFGS}.
The code used to generate the examples shown here is available at
\begin{center}
\url{https://github.com/agussomacal/ConDiPINN}
\end{center}
The interested reader can reproduce our results and test the impact of the variations of certain parameters such as $\e$, $K$, the sampling method, and the machine precision. 

\subsection{Discussion}
\subsubsection{Impact of the number $K$ of training points}
 In this section we discuss the impact of the number $K$ of training points. We fix the machine precision to Float32, and the uniform sampling $-u$.
 
\begin{figure}[ht!]
\centering
\includegraphics[width=0.45\textwidth]{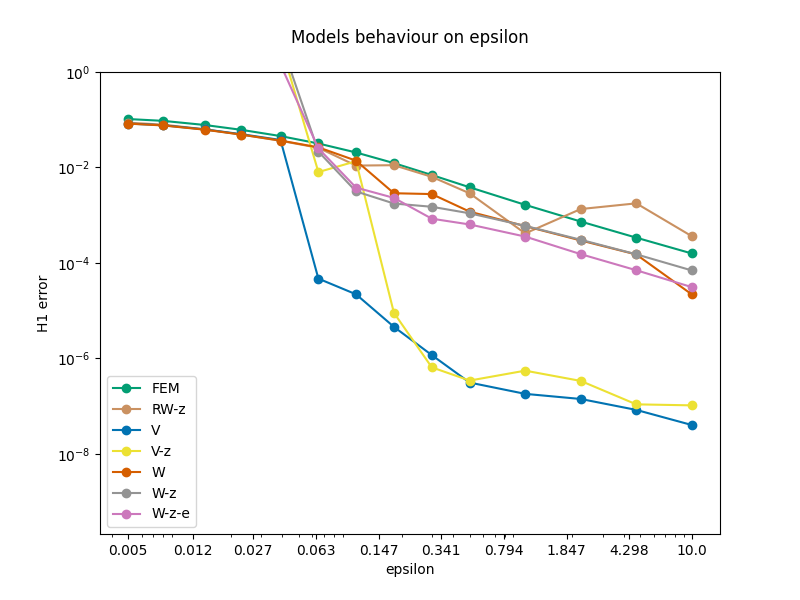}
\includegraphics[width=0.45\textwidth]{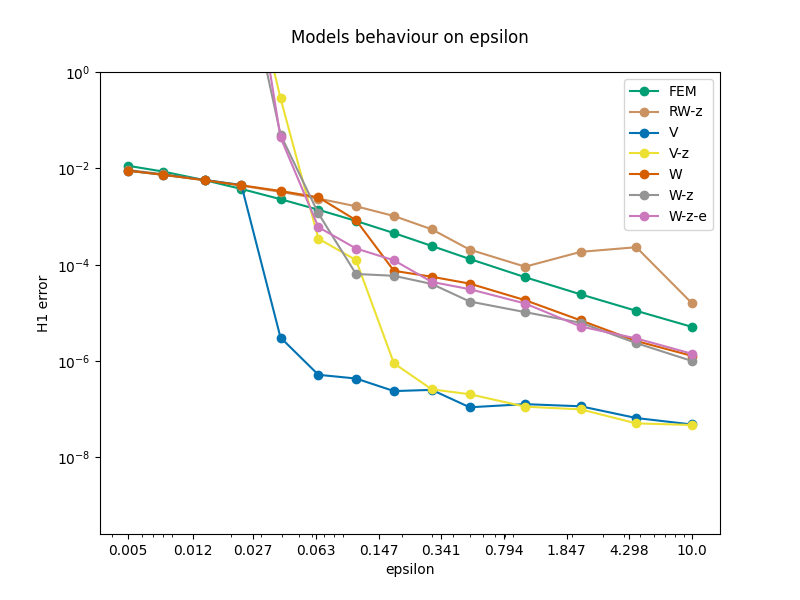}
\includegraphics[width=0.45\textwidth]{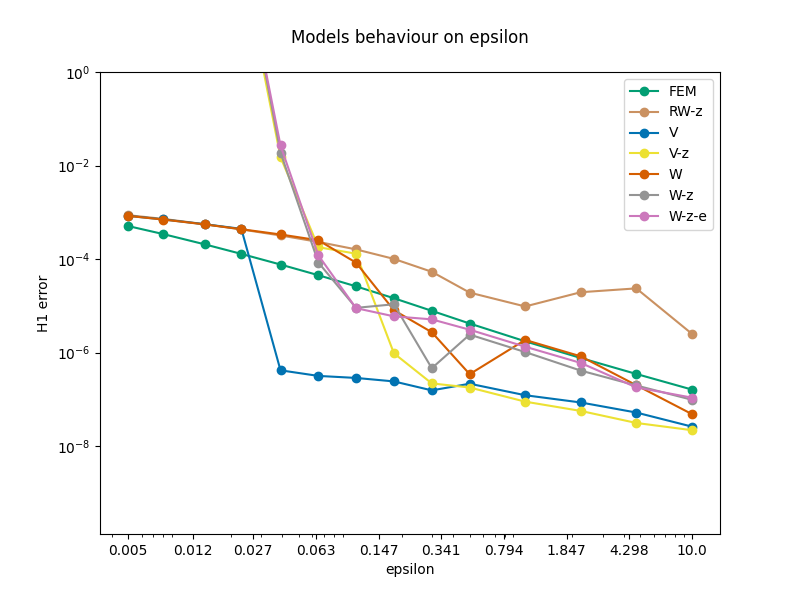}
\includegraphics[width=0.45\textwidth]{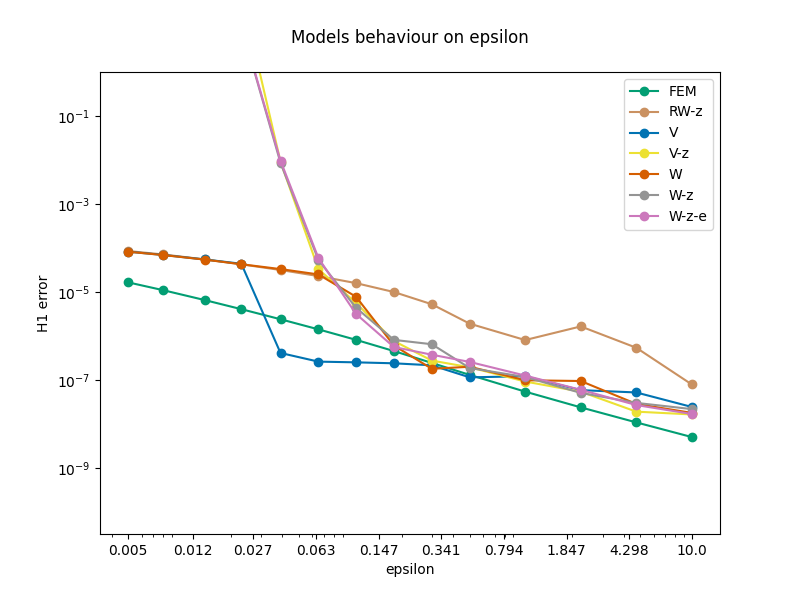}
\caption{Comparison of the behavior of the $h^1$ error for the different methods and different number of sampling points in training. From top to bottom and from left to right, the first figure is produced for $K=10$, the second for $K=100$, the third for $K=1000$ and the last one for $K=10000$. The set of points to train and test have been chosen with the \textit{uniform} sampling method. The precision has been fixed to Float32 for all the tests.}
\label{fig:linspace_diffNtrains}
\end{figure}

\begin{figure}[ht!]
\centering
\includegraphics[width=0.45\textwidth]{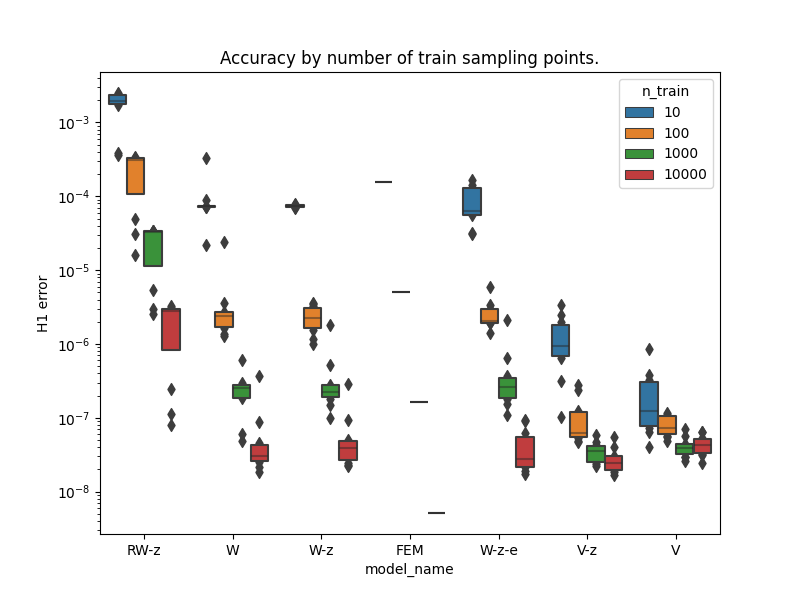}
\includegraphics[width=0.45\textwidth]{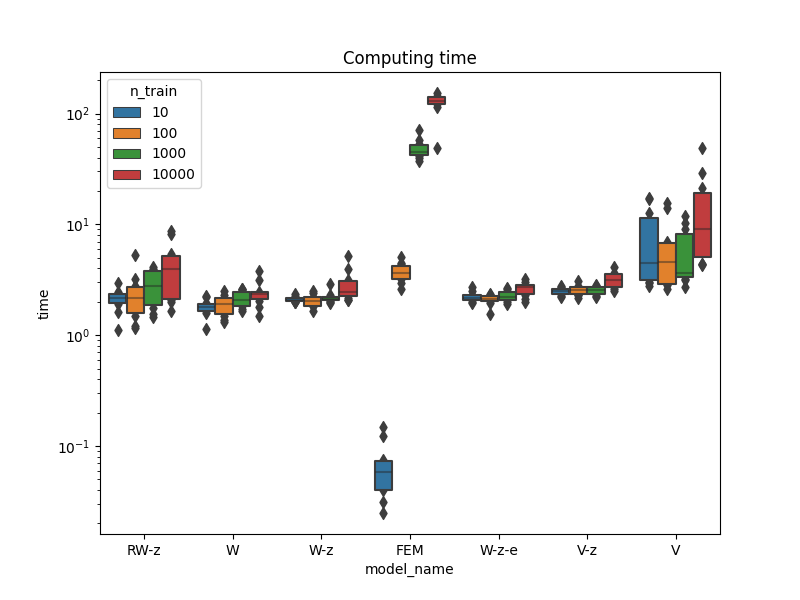}
\caption{For $\e=10$ (region where all methods work well), we look at the comparison between methods and the difference with respect to the number of training points K. The $h^1$ error (right) and the computation times (left). The set of points to train and test have been chosen with the \textit{uniform} sampling method. The precision has been chosen as Float32 for all the tests.}
\label{fig:linspace_diffNtrains_times}
\end{figure}

 Figure \ref{fig:linspace_diffNtrains} shows the best result obtained in the tests, i.e., the minimum value of the $h^1$ norm obtained in the 10 different simulations, plotted against the values of $\e$. In Figure \ref{fig:linspace_diffNtrains_times}, we fix $\epsilon=10$, and plot statistics on the accuracy $e_{h^1}$ (left plot) and computation runtimes for different $K$ (right plot), and for the different methods.
 
 From these figures, we first notice that the approximation of FEM degrades when $\e$ decreases. However, the quality globally improves when the number of training point increases (see, e.g., Figure \ref{fig:linspace_diffNtrains} - left plot). The rate of improvement is  quadratic as we can see from the right plot in Figure \ref{fig:linspace_diffNtrains_times}. In addition, when looking at the runtimes (Figure \ref{fig:linspace_diffNtrains_times} - right) we observe the expected quadratic increase with respect to the number $K$ of training points.
 
 We can next study the behavior of Vanilla PINN and compare to FEM. We observe that it performs at around constant accuracy for any number of training points until around $\e=0.027$ where stops producing reliable approximations (see Figure \ref{fig:linspace_diffNtrains}). One remarkable observation is that the Vanilla PINN error for large values of $\e$ and small number of training points $K=10$ is comparable to the FEM errors with a much larger number of degrees of freedom $K>10^3$ (see left plot in Figure \ref{fig:linspace_diffNtrains_times}). Regarding the runtime to fit the neural network, we see that it is roughly constant for all values of $K$ and it is comparable to the runtime of the FEM method with $K=100$ (right plot in Figure \ref{fig:linspace_diffNtrains_times}).
 
 We next comment on the other PINN-based variational methods. For $\e$ large enough, we observe that all the variational based methods follow the same error trend as FEM both with respect to $\e$ and $K$ and for $K<10^4$ they even perform marginally better. With respect to the computing time, all the methods perform with almost constant time with respect to $K$ and similarly to a FEM method with $K=100$ degrees of freedom. However, for $\e < 0.63$, the methods $W-z$, $W-z-e$ and $V-z$ blow up and lose completely their approximation capabilities. We conjecture that this is due to the fact that the neural network is used to approximate the solution $z$ from the transformed problem, and there is an exponential term to go back from $z$ to $u$
 (see equation \eqref{eq:change-var}). This may lead to machine precision overflows (in the exponential computation) and underflows (the neural network has to learn very small values of $z$ which also are in the limits of precision). To address this issue, we have explored two possible strategies: one was by directly minimizing over $u$ while maintaining the weak formulation which accounts for the method $W$. The second approach is to perform the re-scaling of the domain $RW-z$. In both cases the blow up caused by the exponential is solved although the re-scaling method $RW-z$ doesn't perform as good as others in the region with large $\e$ values.
 
\begin{figure}[ht!]
\centering
\includegraphics[width=0.32\textwidth]{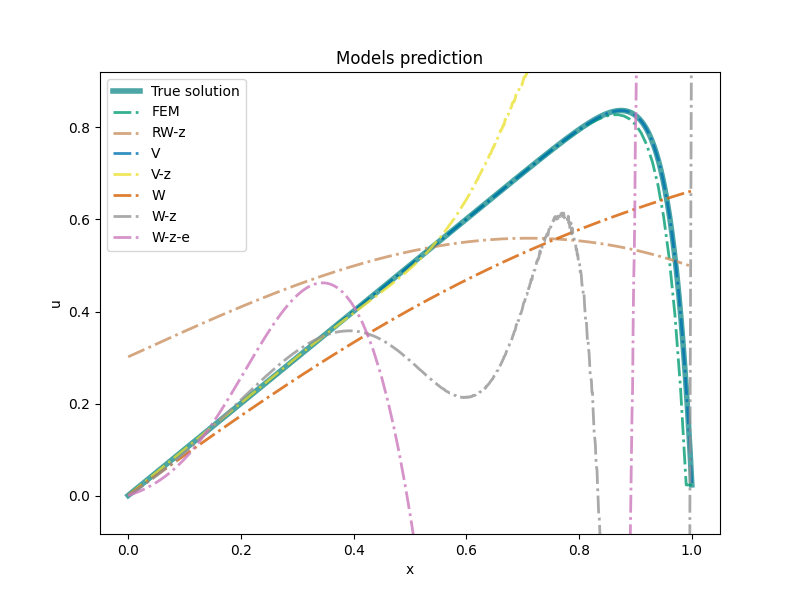}
\includegraphics[width=0.32\textwidth]{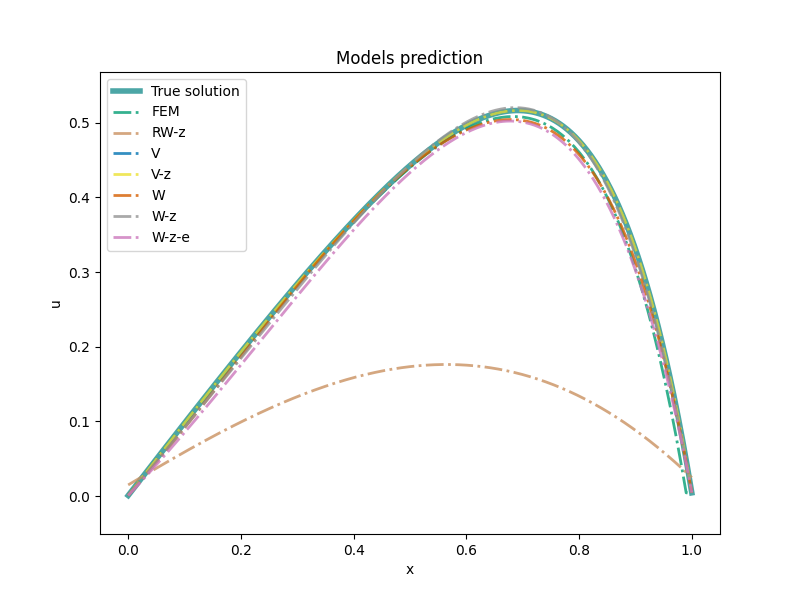}
\includegraphics[width=0.32\textwidth]{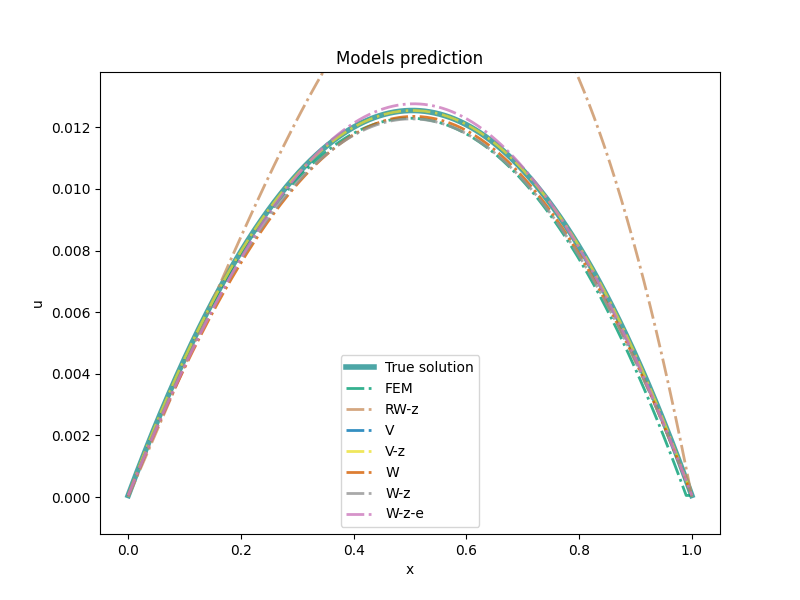}
\caption{The best approximated solution out of $10$ repetitions, for each model, and with $K=100$ training samples. From left to right: $\varepsilon=0.039$, $0.18$, $10$. The interest reader may experiment other configurations in our provided code.}
\label{fig:predictions}
\end{figure}

 We finish this section by plotting in Figure \ref{fig:predictions} the best approximated solution for each model, and different values of $\e$. The interested reader may experiment other configurations in our provided code. The most striking observation is that only FEM and the vanilla PINN method recover the final shape of the exact solution when $\e$ is small. The other variational PINN methods fail despite that some of them exhibit comparable values to FEM in the generalization errors as Figure \ref{fig:linspace_diffNtrains} illustrates. This observations suggests that perhaps other types of error metrics should be introduced in order to be able to better distinguish between ``good solution shapes'' and ``bad ones''. 
 
\subsubsection{Impact of Machine Precision}
\label{sec:FloatPrecision}
\begin{figure}[ht!]
\centering
\includegraphics[width=0.4\textwidth]{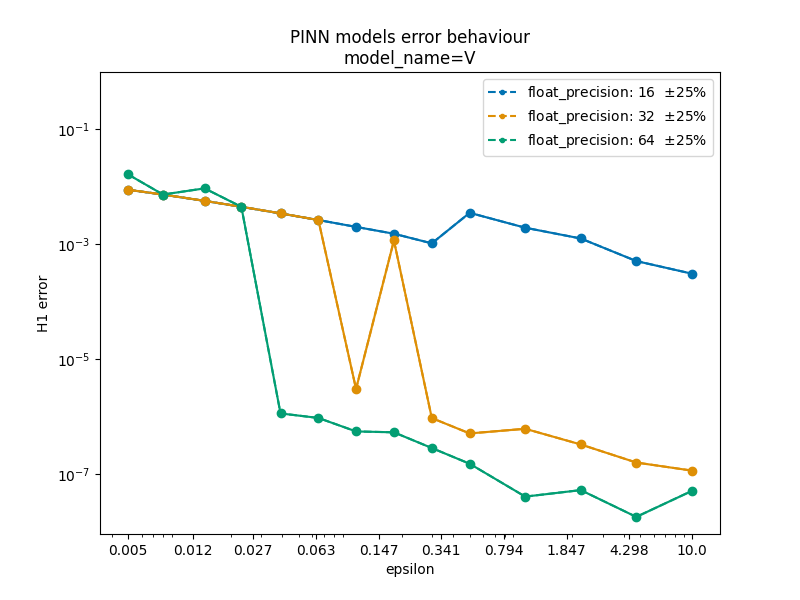}
\includegraphics[width=0.4\textwidth]{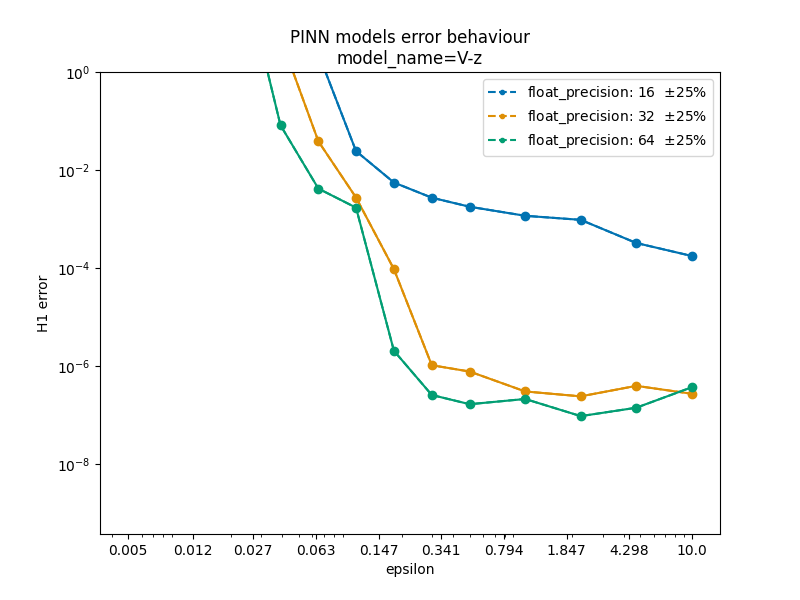}
\includegraphics[width=0.4\textwidth]{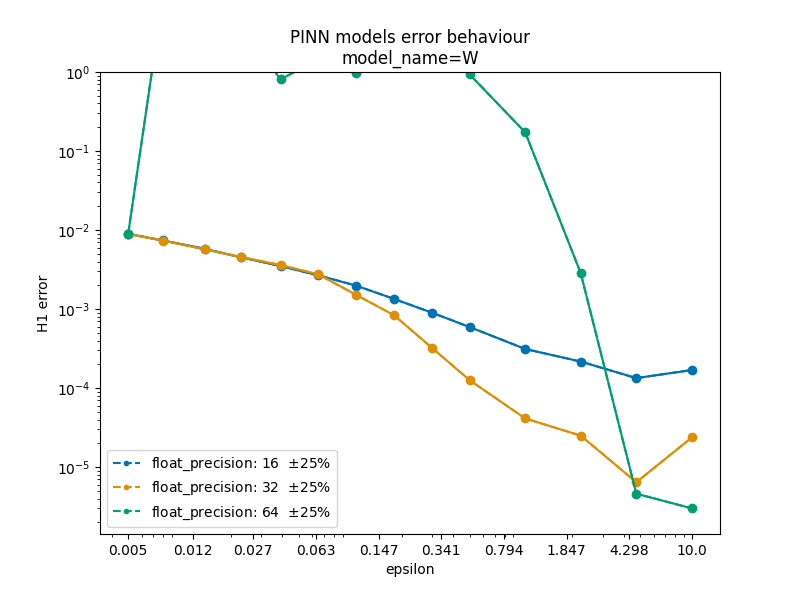}
\includegraphics[width=0.4\textwidth]{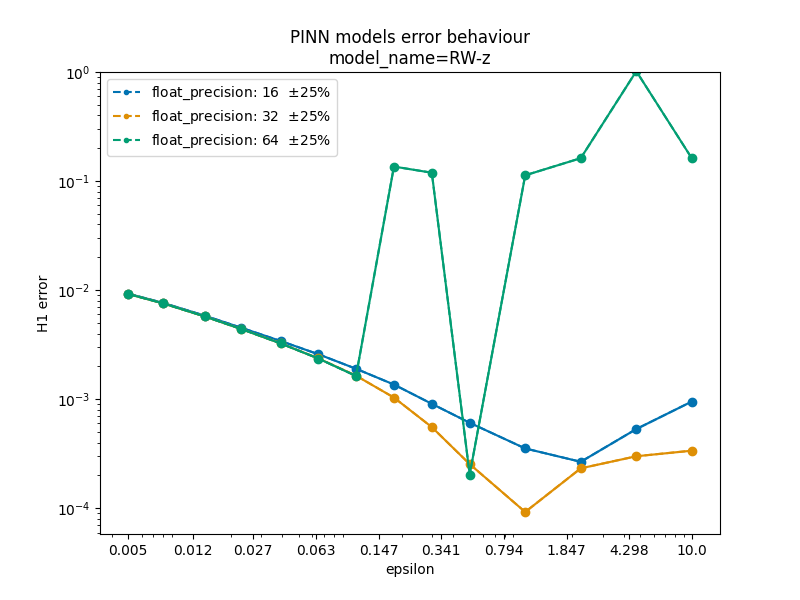}
\includegraphics[width=0.4\textwidth]{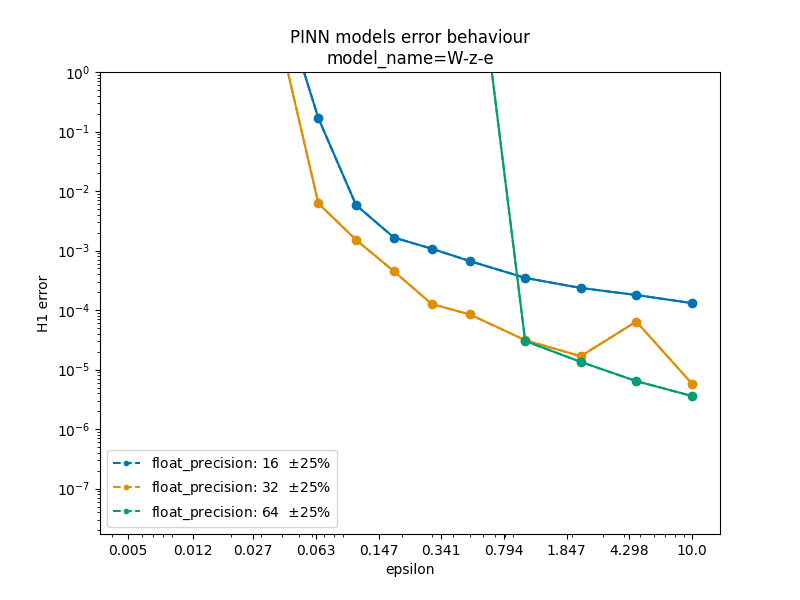}
\includegraphics[width=0.4\textwidth]{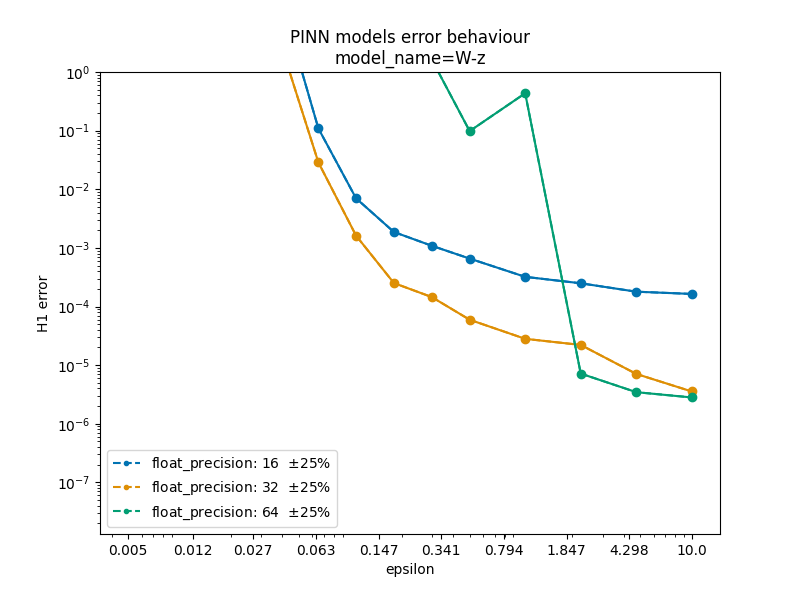}
\caption{Here, the comparison of the behavior of the model for different float precision. The tests have been performed for $K=100$ and uniform sampling.}
\label{fig:floatprecision}
\end{figure}

Figure \ref{fig:floatprecision} shows the $h^1$-error of the different approximated solution by changing the  machine precision in the parameters of the neural networks for the different values of $\e$: Float16, Float32 and Float64. There is an improvement when going from Float16 to Float32 in all methods. Interestingly, we did not obtain very satisfactory results when working with Float64 precision. This precision seems to difficult the convergence to good quality minima: even after 10 repetitions, we failed to find good results. However, as the plots show, when a good minimum is found, it delivers slightly better approximation than lower machine precisions. For these reasons we have performed our experiments using the Float32 which seemed the most stable choice.

\subsubsection{Impact of Sampling Strategy}
\label{sec:SamplingStrategy}
\begin{figure}[ht!]
\centering
\includegraphics[width=0.45\textwidth]{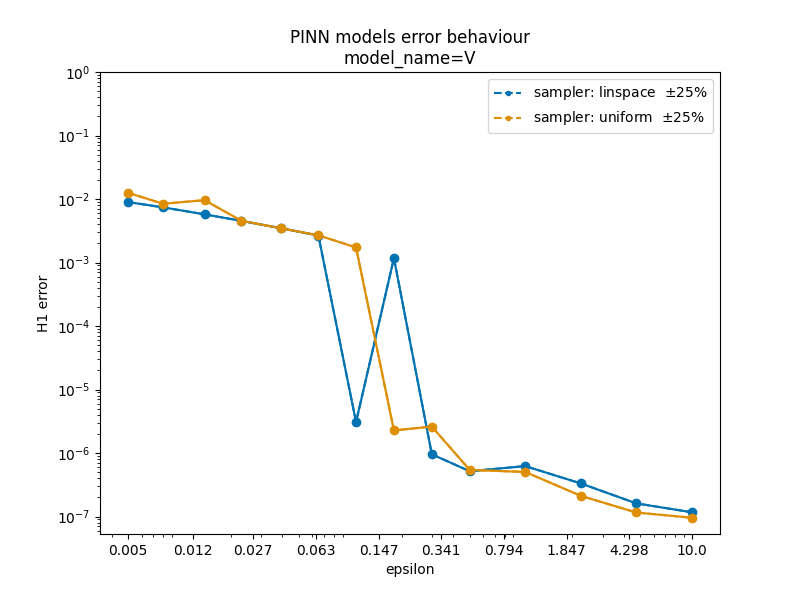}
\includegraphics[width=0.45\textwidth]{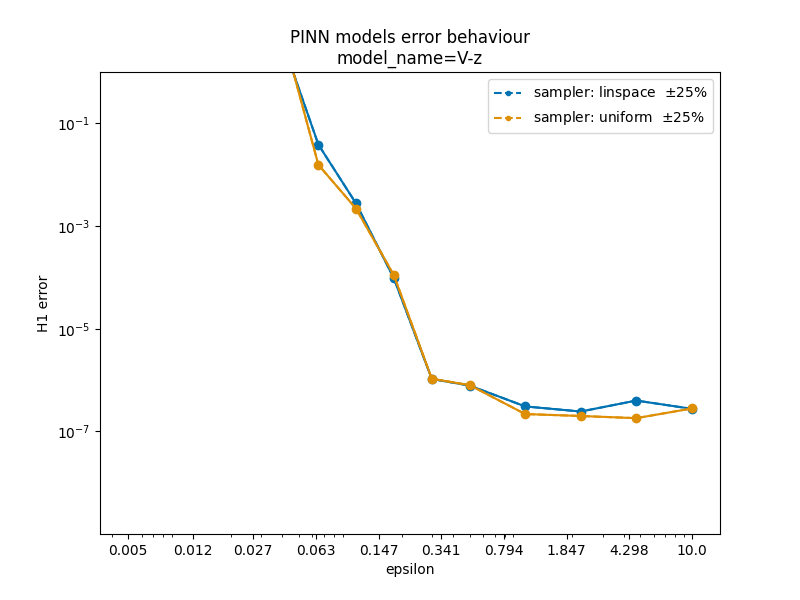}
\includegraphics[width=0.45\textwidth]{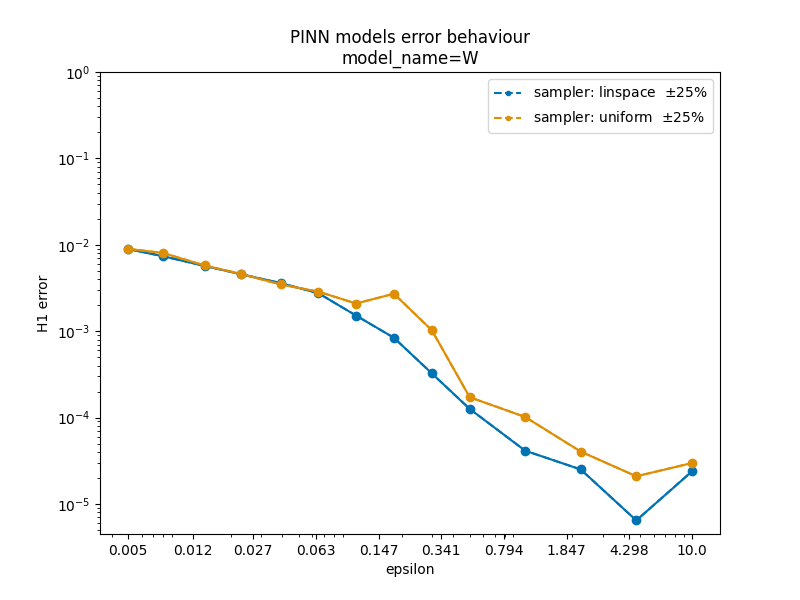}
\includegraphics[width=0.45\textwidth]{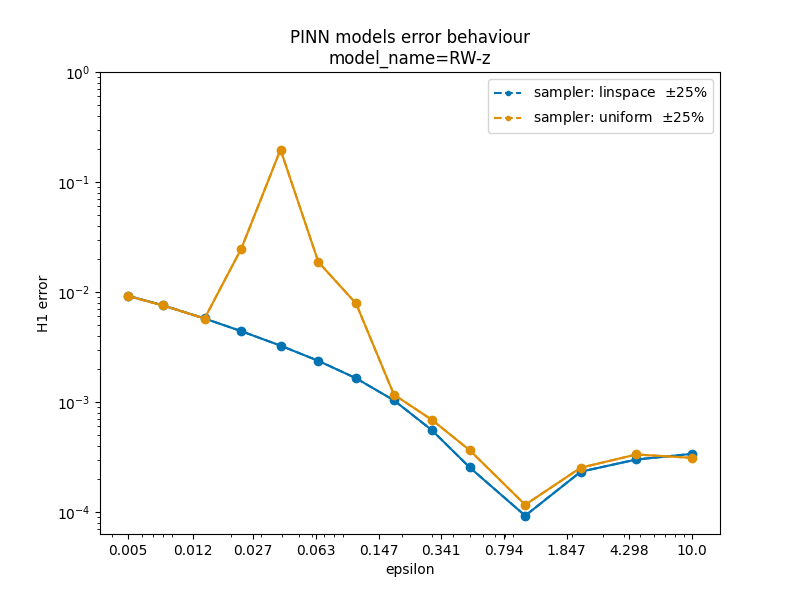}
\includegraphics[width=0.45\textwidth]{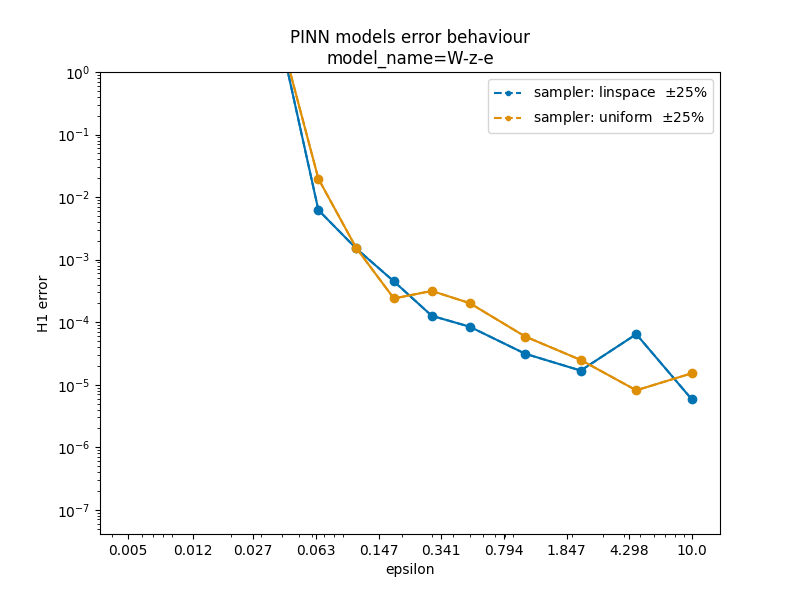}
\includegraphics[width=0.45\textwidth]{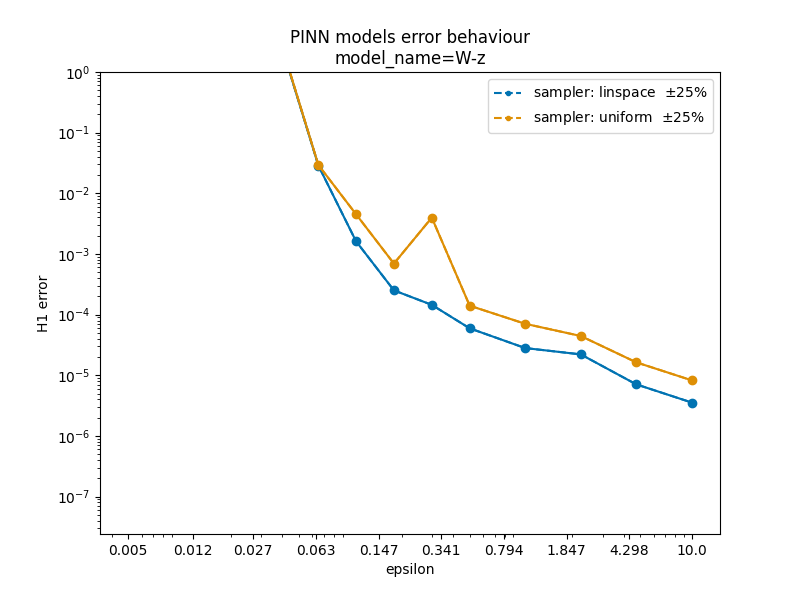}
\caption{Here, the comparison of the behavior of the model for the two different sampling strategies. The tests have been performed for $K=100$ and the float precision equal to Float32.}
\label{fig:sampling_strategy}
\end{figure}

Figure \ref{fig:sampling_strategy} shows the $h^1$-error of the different approximated solution by changing the sampling strategy. For all models, the $uniform$ strategy is found to be either as good as the $random$ or slightly better. For this reason we performed all the experiments using the $uniform$ strategy.

\subsection{Conclusions from the numerical experiments}

The above numerical experiments depict a contrasted landscape concerning the merits and limitations of deep learning-based approaches when the solutions become low regular:

\begin{itemize}
\item For large values of $\e$ when solutions are rather regular, some PINNs perform clearly better than FEM regarding the generalization errors. The superiority is particularly remarkable for very small number $K$ of training points. However, the shapes of PINN solutions are sometimes not as satisfactory as the ones given by FEM.
\item  For the challenging case where $\e$ becomes small and solutions become less regular (which was the main motivation of our study), the accuracy of the variational neural-network methods is essentially comparable or worse to the one given by FEM in terms of generalization errors. Some PINN variational approaches become too unstable and the errors blow up. Only FEM and the vanilla PINN approach seem to be able to recover the correct shape of the exact function. The latter one has however the risk of sometimes falling into local minima with bad shapes. 
\item The runtimes are clearly in favor to PINN methods as Figure \ref{fig:linspace_diffNtrains_times} illustrates, and the simplicity of implementation is also in favor to all PINN methods.
\end{itemize}

\section{Future research directions and extensions}
One important point to explore in future works concerns the choice of the loss function for the training, and also the metric to evaluate generalization errors. It will also be interesting to explore if adaptive sampling strategies during the training could help to recover good solutions in a more stable manner. Finally, the impact of the machine precision in some steps involving exponential transformations seems also to be an important obstacle to retrieving stable solutions. It would be interesting to develop strategies that circumvent this issue. All these developments will play a crucial role in order to address higher dimensional problems with similar characteristics as the one considered here.

\newpage
\appendix

\section{Analytic solution in dimension 1}
\label{sec:AppendixAnaliticSolution}

The aim of this section is to give the analytic expression of the solution of (\ref{general_equation})-(\ref{Robin}) in the case when $d=1$, $\Omega = (0, 1)$, and $F$ and $f$ are assumed to be equal to some constant real numbers. Thus, in this section, using a slight abuse of notation, we assume that $F,f,g \in \mathbb{R}$. Let us also introduce $g_0, g_1\in \mathbb{R}$ so that $g(0) = g_0$ and $g(1) = g_1$. The problem then reads as follows: find $u: (0,1)\to \mathbb{R}$ solution to
\begin{equation}
\left\{
\begin{array}{rl}
- \e u''(x) + F u'(x) = f, & \quad \forall x\in (0, 1),\\
 - \alpha u'(0) + \kappa u(0) = g_0, & \\
 \alpha u'(1) + \kappa u(1) = g_1. &
\end{array}\right. 
\end{equation}
Then, it can be easily checked that the solution to this equation reads as
$$
u(x) = C_1 + C_2 e^{\frac{Fx}{\e}} + \frac{f}{F} x
$$
where $C_1$ and $C_2$ are constants that are determined with the Robin boundary conditions. They satisfy the system
\begin{equation*}
\begin{pmatrix}
\kappa & \kappa - \frac{\alpha F}{\e}  \\
\kappa & \kappa e^{\frac{F}{\e}} + \alpha \frac{F}{\e} e^{\frac{F}{\e}}
\end{pmatrix}
\begin{pmatrix}
C_1 \\
C_2
\end{pmatrix}
=
\begin{pmatrix}
g_0 + \alpha \frac{f}{F} \\
g_1 - \frac{f}{F} ( \kappa+ \alpha)
\end{pmatrix}
\end{equation*}
which invertible except for
$$
\kappa = 0, \quad \text{or} \quad
\frac{\alpha}{\kappa} = \frac{\e (1 - e^{\frac{F}{\e}})}{F(1+ e^{\frac{F}{\e}})}.
$$ 
In the following, we assume that the values of $\kappa$ and $\alpha$ do not take these values, and the above system is invertible.

\newpage

\section{l2 error plots}
\label{sec:l2error}

\begin{figure}[ht!]
\centering
\includegraphics[width=0.45\textwidth]{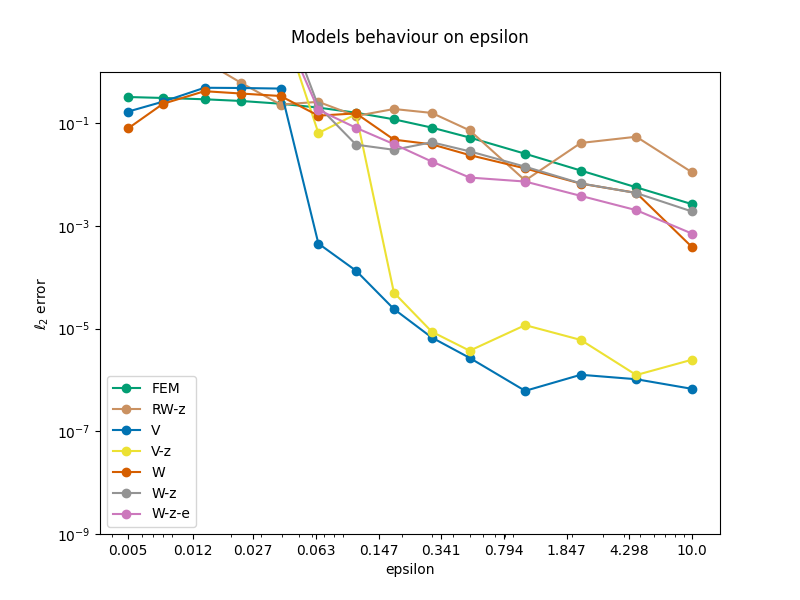}
\includegraphics[width=0.45\textwidth]{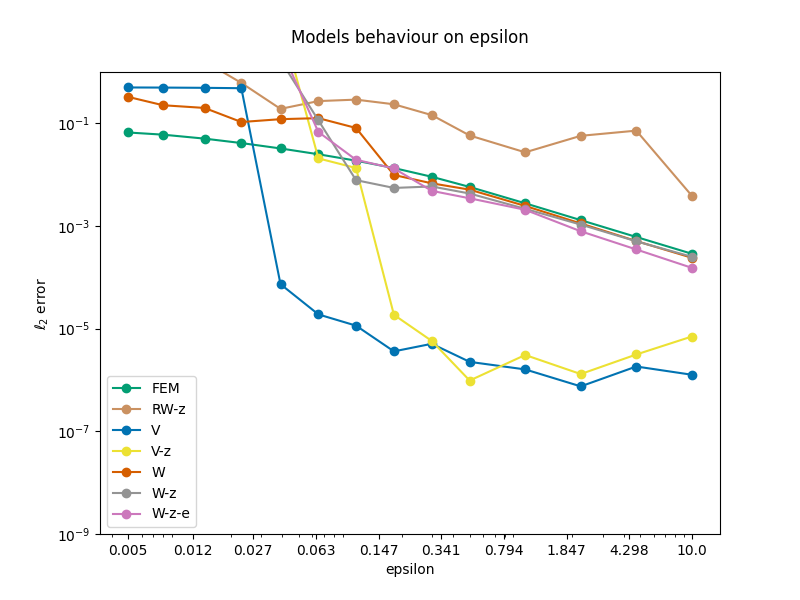}
\includegraphics[width=0.45\textwidth]{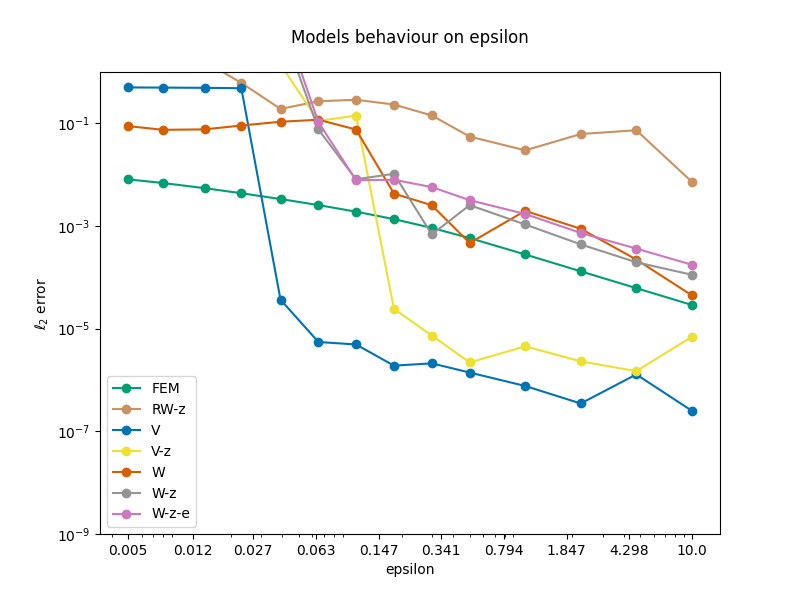}
\includegraphics[width=0.45\textwidth]{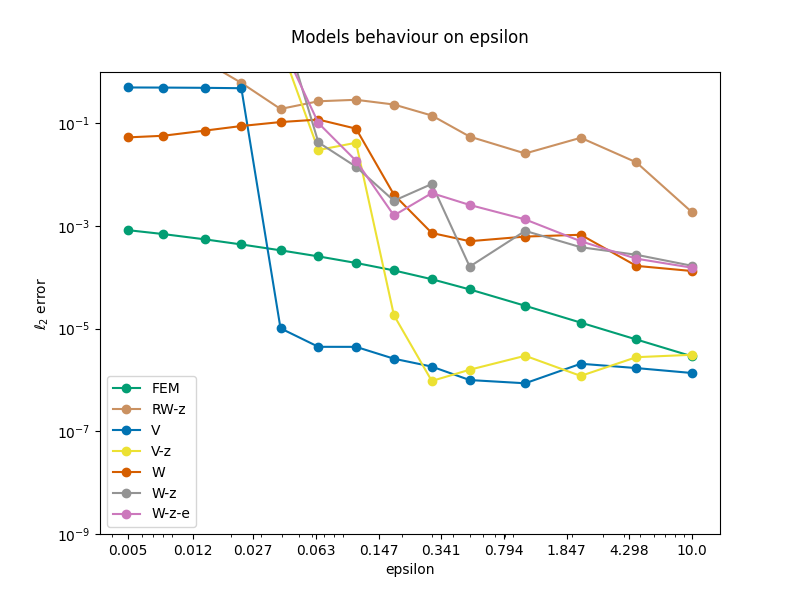}
\caption{The comparison of the behavior of the $l^2$ error for the different methods and different number of sampling points in training. From up to down and from left to right, the first figure is produced for $K=10$, the second for $K=100$, the third for $K=1000$ and the last one for $K=10000$. The set of points to train and test have been chosen with the \textit{uniform} sampling method. The precision has been chosen as Float32 for all the tests.}

\end{figure}

\newpage

\bibliographystyle{unsrt}
\bibliography{literature}

\end{document}